\documentclass[11pt]{article}
\usepackage{geometry}                
\usepackage{graphicx}
\usepackage{amssymb}
\usepackage{epstopdf}
\usepackage{amsmath,amsthm,amscd,amssymb}
\usepackage{latexsym}
\numberwithin{equation}{section}

\theoremstyle{plain}
\newtheorem{theorem}{Theorem}[section]
\newtheorem{lemma}[theorem]{Lemma}
\newtheorem{corollary}[theorem]{Corollary}

\newtheorem{conjecture}[theorem]{Conjecture}

\theoremstyle{definition}
\newtheorem{definition}[theorem]{Definition}

\theoremstyle{remark}
\newtheorem{remark}[theorem]{Remark}

\newtheorem{case[theorem]}{Case}

\title{Erdos distance problem in vector spaces over finite fields}
\author{Alex Iosevich and Misha Rudnev}

\begin{document}
\maketitle

\begin{abstract} We study the Erd\"os/Falconer distance problem in vector spaces over finite fields. Let ${\Bbb F}_q$ be a finite field with $q$ elements and take $E \subset {\Bbb F}^d_q$, $d \ge 2$. We develop a Fourier analytic machinery, analogous to that developed by Mattila in the continuous case, for the study of distance sets in ${\Bbb F}^d_q$ to provide estimates for minimum cardinality of the distance set $\Delta(E)$ in terms of the cardinality of $E$. Kloosterman sums play an important role in the proof. \end{abstract}

\tableofcontents

\newpage

\section{Introduction}

In the Euclidean setting, the Erdos distance conjecture says that if $E$ is a finite subset of ${\Bbb R}^d$, $d \ge 2$, then 
\begin{equation} \label{euclideanconjecture} \# \Delta_{{\Bbb R}^d}(E) \gtrapprox {(\# E)}^{\frac{d}{2}}, \end{equation} where 
$$ \Delta_{{\Bbb R}^d}=\{|x-y|: x,y \in E\}, $$ and 
$$ {|x-y|}^2={(x_1-y_1)}^2+\dots+{(x_d-y_d)}^2.$$ 

Here, and throughout the paper, $X \lesssim Y$ means that there exists $C>0$ such that $X \leq CY$, $X\gtrsim Y$ means $Y\lesssim X$, and $X\approx Y$ if both $X\lesssim Y$ and $X\gtrsim Y$. Besides, $X \lessapprox Y$ means that for every $\epsilon>0$ there exists $C_{\epsilon}>0$ such that $X \leq C_{\epsilon}q^{\epsilon}Y$, where $q$ is a large controlling parameter. 

See, for example, \cite{Ma02} for the description of the Erdos distance problem in Euclidean space and references to recent results. We mention in passing that the Erdos distance conjecture is not solved in any dimension, in Euclidean or any other setting. The best known result in the Euclidean plane is due to Katz and Tardos (\cite{KT04}) who prove that 
$$\# \Delta_{{\Bbb R}^2}(E) \gtrsim {(\# E)}^{\approx .86}.$$ 

In this paper we study the Erdos distance problem in vector spaces over finite fields. Let ${\Bbb F}_q$ denote the finite field with $q$ elements, and let ${\Bbb F}^d_q$ denote the $d$-dimensional vector space over this field. Let $E \subset {\Bbb F}_q^d$, $d \ge 2$. Then the analog of the classical Erdos distance problem is to determine the smallest possible cardinality of the set 
$$ \Delta(E)=\{{{|x-y|}^2=(x_1-y_1)}^2+\dots+{(x_d-y_d)}^2: x,y \in E\}, $$ viewed as a subset of ${\Bbb F}_q$. 

In the finite field setting, the estimate (\ref{euclideanconjecture}) cannot hold without further restrictions. To see this, let $E={\Bbb F}^d_q$. Then $\# E=q^d$ and $\# \Delta(E)=q$. With this example as our guide, we are led to the following conjecture. 

\begin{conjecture} \label{conjecture} Let $E \subset {\Bbb F}^d_q$ of cardinality $\lesssim q^{\frac{d}{2}}$. Then 
$$ \# \Delta(E) \gtrsim {(\# E)}^{\frac{2}{d}}. $$ \end{conjecture} 

A Eucludean plane argument due to Erd\"os (\cite{E45}) can be applied to the final field set-up to show that if $d=2$, then
\begin{equation} \label{trivial2d} \# \Delta(E) \gtrsim {(\# E)}^{\frac{1}{2}}. \end{equation} 

In higher dimensions matters are more subtle in the context of vector spaces over finite fields because intersection of analogs of spheres in ${\Bbb F}^d_q$ may be quite complicated, and the standard dimensional induction in ${\Bbb R}^d$ argument (see e.g. \cite{AP95}) that allows one to bootstrap the estimate (\ref{trivial2d}) into the estimate 
\begin{equation} \label{trivialhighd} \# \Delta_{{\Bbb R}^d}(E) \gtrsim {(\# E)}^{\frac{1}{d}} \end{equation} does not immediately go through. We establish the finite field analog of the estimate (\ref{trivialhighd}) below using Fourier analytic methods and number theoretic properties of Kloosterman sums. 

The idea here is that while the Euclidean form for the Erd\"os distance conjecture has no chance of holding for certain large subsets of ${\Bbb F}^d_q$, it may in fact hold for subsets of smaller size. A similar idea is explored in \cite{BKT04} in the context of incidence theorems. 

Another way of thinking of Conjecture \ref{conjecture} is in terms of the Falconer distance conjecture (\cite{Fal86}) in the Euclidean setting which says that if the Hausdorff dimension of a set in ${\Bbb R}^d$ exceeds $\frac{d}{2}$, then the Lebesgue measure of the distance set is positive. Conjecture \ref{conjecture} implies that if the size of the set is greater than $q^{\frac{d}{2}}$, then the distance set contains a positive proportion of all the possible distances, an analogous statement. 

The methods of this paper are strongly motivated by the Falconer conjecture. In particular, a significant part of this paper is dedicated to the derivation of the finite field analog of Fourier theory for distance sets initially developed in the continuous setting by Falconer (\cite{Fal86}) and Mattila (\cite{Mat87}). See also some recent progress on this problem due to Bourgain (\cite{B94}), Erdogan (\cite{Erd05}) and Wolff (\cite{W99}). The best currently known result is due to Erdogan (\cite{Erd05}) and Wolff (\cite{W99}) who proved that the Lebesgue measure of the distance set is positive provided that the Hausdorff dimension of the set exceeds $\frac{d}{2}+\frac{1}{3}$. Note that Theorem \ref{trivial} below corresponds to the exponent $\frac{d}{2}+\frac{1}{2}$, proved in the continuous case by Falconer (\cite{Fal86}). The proof in the finite field case is more difficult and involves non-trivial number theory, mainly hidden in the known estimates for Kloosterman sums. 

\subsection{Statement of results} 

\begin{definition} The Fourier transform of a function $F: {\Bbb F}^d_q \to {\Bbb F}_q$ is given by 
$$ \widehat{F}(m)=q^{-d} \sum_{x \in {\Bbb F}^d_q} e^{-\frac{2 \pi i x \cdot m}{q}} F(x), $$ for $m \in {\Bbb F}^d_q$, \end{definition} where ${\Bbb F}_q$ is identified with the roots of unity on the circle in the usual way. See, for example, \cite{G03} and \cite{StSh03}. 

See the following section for the description of basic properties of the Fourier transform in vector spaces over finite fields. 

\begin{theorem} \label{mattila} Let $E \subset {\Bbb F}_q^d$, $d \ge 2$. Suppose that $\# E \gtrsim q^{\frac{d}{2}}$. Let 
\begin{equation} \label{mattilaintegral} {\cal M}(q)=\frac{q^{3d+1}}{{(\# E)}^4} \sum_{\{(m,m') \in {\Bbb F}^d_q \times {\Bbb F}^d_q: {|m|}^2={|m'|}^2\}} {|\widehat{E}(m)|}^2 {|\widehat{E}(m')|}^2. \end{equation}

Then 
$$ \# \Delta(E) \gtrsim \min \left\{q, \frac{q}{{\cal M}(q)} \right\}.$$ 
\end{theorem} 

See Section 6 below for examples of the quantity ${\cal M}(q)$ computed for various natural subsets of ${\Bbb F}^d_q$ of critical cardinality. 

\begin{remark} The quantity ${\cal M}(q)$ is the finite field analog of the Mattila integral (see \cite{Mat87}), given by 
\begin{equation} \label{contmat} \int_1^{\infty} {\left( \int_{S^{d-1}} {|\widehat{\mu}(t \omega)|}^2 d\omega \right)}^2 t^{d-1}dt, \end{equation} where $\mu$ is a Borel measure on a set $E \subset {\Bbb R}^d$, $d \ge 2$, of Hausdorff dimension $>\frac{d}{2}$. Mattila proves that if this quantity is bounded for some Borel measure $\mu$ supported on $E$, then the distance set of $E$ has positive Lebesgue measure. \end{remark} 

\begin{definition} In analogy with the Euclidean case, we say that $E \subset {\Bbb F}^d_q$ is a Salem set if for every non-zero element $m$ of ${\Bbb F}^d_q$, 
\begin{equation} \label{salemest} |\widehat{E}(m)| \lesssim q^{-d} \cdot \sqrt{\# E}. \end{equation} \end{definition} 

See Lemma \ref{salemexample} and Lemma \ref{salemsphere} for some natural examples of Salem sets. In particular, we shall see that the "paraboloid" 
$$ P=\{(x, {|x|}^2) \in {\Bbb F}^{d-1}_q \times {\Bbb F}_q\},$$ and the "sphere" 
$$ S_r=\{x \in {\Bbb F}^d_q: {|x|}^2=r\}$$ are Salem sets. We note in passing that continuous analogs of these objects are Salem sets in the Euclidean setting. 

\begin{theorem} \label{salem} Suppose that $E \subset {\Bbb F}^d_q$ is a Salem set of cardinality $\gtrsim q^{\frac{d}{2}}$. Then the Conjecture \ref{conjecture} holds. \end{theorem}

\begin{theorem} \label{trivial} Let $E \subset {\Bbb F}^d_q$. Then 
$$  \# \Delta(E) \gtrsim \min \left\{q, \frac{\# E}{q^{\frac{d-1}{2}}} \right\}. $$
\end{theorem}

\begin{corollary} \label{triviality} Suppose that $\# E \gtrsim q^{\frac{d}{2}}$. Then the estimate (\ref{trivialhighd}) holds. \end{corollary} 

\begin{corollary} \label{lesstrivial} Suppose that $\# E \approx q^{\frac{d+1}{2}}$. Then $\# \Delta(E) \gtrsim {(\# E)}^{\frac{2}{d+1}}.$ \end{corollary}

\begin{remark} Kloosterman sums play an important role in the proof of Theorem \ref{trivial}. We also establish the fact that the sphere is a Salem set in the course of the argument. \end{remark}

The bulk of the work on the continuous analog of the problems has centered around the point-wise estimation of the quantity 
\begin{equation} \label{contav} \int_{S^{d-1}} {|\widehat{\mu}(t\omega)|}^2 d\omega, \end{equation} and then plugging the result into the quantity (\ref{contmat}). The discrete analog of (\ref{contav}) is the quantity 
\begin{equation} \label{disav} \sigma^2_E(m)=\sum_{\{m' \in {\Bbb F}^d: {|m'|}^2={|m|}^2\}} {|\widehat{E}(m')|}^2. \end{equation}

\begin{theorem} \label{plug} Suppose that 
\begin{equation} \label{disavest} |\sigma^2_E(m)| \lesssim q^{-\beta}. \end{equation} 

Then 
$$ \# \Delta(E) \gtrsim \min \left\{q, \frac{{(\# E)}^3}{q^{2d-\beta}} \right\}.$$ 
\end{theorem} 

The following is an immediate consequence of Theorem \ref{plug}. 

\begin{corollary} \label{pieinthesky} Suppose that $\# E \gtrsim q^{\frac{d}{2}}$ and the estimate (\ref{disavest}) holds with $\beta=\frac{d}{2}+1$. Then $\# \Delta(E) \gtrsim q$. \end{corollary}

We also have the following positive result. 

\begin{theorem} \label{something} We have 
$$ \sigma^2_E(m) \lesssim q^{-\frac{d+1}{2}} \frac{{(\# E)}^2}{q^d}. $$
\end{theorem} 

The following analog of Falconer's (Euclidean) theorem follows by combining Theorem \ref{plug} and Theorem \ref{something}. 

\begin{corollary} \label{falconer} Suppose that $\# E \gtrsim q^{\frac{d+1}{2}}$. Then $\# \Delta(E) \gtrsim q$. \end{corollary}

\subsection{Related work} Finite field analogs of various theorems in harmonic analysis and geometric combinatorics have been explored in a number of recent papers. See \cite{MT04}  for a finite field version of the restriction phenomenon, and \cite{BKT04} for the discussion of incidence theorems this setting. Also see the references contained in these papers to related work in additive number theory. 

\vskip.125in

\section{Finite field analog of the Fourier transform and applications to distance sets} 

We start out with a quick review of basic definitions and results about the Fourier transform in finite fields.  
See \cite{Mat87} for the description of a similar method in the continuous setting. Let $f$ be a function on ${\Bbb F}_q$. Define the $k$th Fourier coefficient of $f$ by the relation
$$ \widehat{f}(k)=\frac{1}{q} \sum_{j=0}^{q-1} e^{-\frac{2 \pi i jk}{q}} f(j). $$

It is not difficult to show that 
$$ f(j)=\sum_{k \in {\Bbb F}_q} \widehat{f}(k) e^{\frac{2 \pi i jk}{q}},$$ and 
\begin{equation} \label{plancherel1d} \sum_{k \in {\Bbb F}_q} {|\widehat{f}(k)|}^2=\frac{1}{q} \sum_{j \in {\Bbb F}_q} {|f(j)|}^2. \end{equation}

Similarly, if $F$ is a function on ${\Bbb F}_q^d$, 
$$ \widehat{F}(m)=\frac{1}{q^d} \sum_{x \in {\Bbb F}^d_q} e^{-\frac{2 \pi i x \cdot m}{q}} F(x), $$
\begin{equation} \label{inversion} F(x)=\sum_{m \in {\Bbb F}^d_q} e^{\frac{2 \pi i x \cdot m}{q}} \widehat{F}(m), \end{equation} and 
\begin{equation} \label{plancherelhigher} \sum_{m \in {\Bbb F}^d_q} {|\widehat{F}(m)|}^2=\frac{1}{q^d} \sum_{x \in {\Bbb F}^d_q} {|F(x)|}^2. \end{equation}

Our approach to the proof of Theorem \ref{mattila} is the finite field variant of Mattila's $L^2$ technique. Define the measure $\nu$ on $\Delta(E)$ by the relation 
\begin{equation} \label{measure} \sum_{j \in {\Bbb F}_q} f(j) \nu(j)=\frac{1}{{(\# E)}^2} \sum _{x,y \in E} f({|x-y|}^2). \end{equation}

To see precisely what this means, we write 
$$  \frac{1}{{(\# E)}^2}  \sum _{x,y \in E} f({|x-y|}^2)$$
$$=\frac{1}{{(\# E)}^2}  \sum_{j \in {\Bbb F}^q_d} \sum_{\{(x,y) \in E \times E: {|x-y|}^2=j\}} f(j). $$ In other words, 
$$ \nu(j)=\frac{1}{{(\# E)}^2}  \# \{(x,y) \in E \times E: {|x-y|}^2=j\},$$ the incidence function which measures how often a single "distance" $j$ occurs. Observe that 
$$ supp(\nu)=\{j: \nu(j) \not=0\}=\Delta(E). $$

Using (\ref{measure}) observe that the Fourier coefficient of $\nu$, 
$$ \widehat{\nu}(k)=q^{-1} \sum_{j \in {\Bbb F}_q} e^{-2 \pi i \frac{jk}{q}} \nu(j)=q^{-1} \frac{1}{{(\# E)}^2}  \sum_{x,y \in E} e^{-\frac{2 \pi i k {|x-y|}^2}{q}}.$$ 

It follows that 
$$ 1={\left(\sum_{j \in {\Bbb F}_q} \nu(j) \right)}^2 \leq \# \Delta(E) \cdot \sum_{j \in {\Bbb F}_q} \nu^2(j)$$ $$=\# \Delta(E) \cdot q \cdot \sum_{k \in {\Bbb F}_q} {|\widehat{\nu}(k)|}^2, $$ so 
$$ \# \Delta(E) \ge \frac{1}{q \cdot \sum_{k \in {\Bbb F}_q} {|\widehat{\nu}(k)|}^2}. $$

Now, 
\begin{equation} \label{distcoef} \widehat{\nu}(k)=q^{-1} \frac{1}{{(\# E)}^2} \sum_{x \in {\Bbb F}^d_q} T^k_qE(x) \cdot E(x), \end{equation} where the operator $T^k_q$ is defined by 
$$ T^k_qh(x)=\sum_{y \in {\Bbb F}^d_q} e^{\frac{2 \pi i k{|x-y|}^2}{q}} h(y), $$ and $E(x)$ is the characteristic function of $E$. 

We are thus led to study the Fourier coefficients 
$$ \widehat{T^k_qE}(m)=q^{-d} \sum_{x,y \in {\Bbb F}^d_q} e^{\frac{2 \pi i(x \cdot m-k{|x-y|}^2)}{q}} E(y), $$ and we must study Gauss sums of the form 
$$ G(m,k)=\sum_{x,y \in {\Bbb F}^d_q} e^{\frac{2 \pi i(x \cdot m-k{|x|}^2)}{q}}, $$ since by a simple change of variables, 
$$ \widehat{T^k_qE}(m)=G(m,k) \widehat{E}(m). $$

\section{Estimation of Gauss sums and examples of Salem sets} 

We have 
$$ \sum_{x_j \in {\Bbb F}_q} e^{\frac{2 \pi i(m_jx_j-kx^2_j)}{q}}$$ $$=e^{\frac{2 \pi i m^2_j}{4kq}} \sum_{x_j \in {\Bbb F}_q} e^{-\frac{2 \pi i k{(x_j-m_j/2k)}^2}{q}}$$
$$=e^{\frac{2 \pi i m^2_j}{4kq}} g(k),$$ where $g(k)$ is the "standard" Gauss sum 
$$ g(k)=\sum_{x_j \in {\Bbb F}_q} e^{\frac{2 \pi i k x^2_j}{q}}. $$

It follows that if $k \not=0$, then 
\begin{equation} \label{Gauss} G(m,k)=e^{\frac{2 \pi i {|m|}^2}{4kq}} g^d(k). \end{equation}

It is well known that 
\begin{equation} \label{gauss} g(k)=\pm i\sqrt{q}, \end{equation} so 
\begin{equation} \label{gaussd} g^d(k)={(\pm i)}^d \cdot q^{\frac{d}{2}}. \end{equation} Indeed, 
$$ {|g(k)|}^2=\sum_{u,v \in {\Bbb F}_q} e^{\frac{2 \pi i k(u^2-v^2)}{q}}$$
$$=\sum_{t \in {\Bbb F}_q} e^{\frac{2 \pi i kt}{q}} n(t),$$ where 
$$ n(t)=\# \{(u,v) \in {\Bbb F}_q \times {\Bbb F}_q: u^2-v^2=t\}.$$ 

\begin{lemma} We have $n(0)=2q-1$, and $n(t)=q-1$ if $t \not=0$. \end{lemma} 

The former is obvious. To see the latter, consider a homomorphism $h: {\Bbb F}_q^{*} \to {\Bbb F}_q^{*}$ given by $h(u)=u^2$, where ${\Bbb F}_q^{*}$ denotes the multiplicative group of ${\Bbb F}_q$. The kernel of $h$ is $\{-1,1\}$. It follows that the image of $h$ has $\frac{q-1}{2}$ elements. In other words, exactly half the elements in ${\Bbb F}_q^{*}$ are squares. This implies the second claim immediately and the proof of the lemma is complete. 

Alternatively, we can write $u^2-v^2=(u-v)(u+v)$. since $u-v$ and $u+v$ determine $u$ and $v$ uniquely, it suffices to count the number of solutions of the equation $u'v'=t$, $t \not=0$. There are $q-1$ choices for $u'$, say, and $v'$ is completely determined. The same outcome as above follows. 

We conclude that 
$$ {|g(k)|}^2=q+(q-1) \sum_{t \in {\Bbb F}_q} e^{\frac{2 \pi i kt}{q}}=q. $$

Suppose that $-1$ is not a square in ${\Bbb F}_q$. It follows that 
$$ g(k)+\overline{g(k)}=\sum_{t \in {\Bbb F}_q} e^{\frac{2 \pi i kt}{q}}+e^{-\frac{2 \pi i kt}{q}}$$ runs over each of the elements of ${\Bbb F}_q$ exactly twice and thus equals $0$. It follows that $g(k)$ is purely imaginary. If $-1$ is a square in ${\Bbb F}_q$, then $\pm i$ is simply replaced by a different constant. See, for example, \cite{Lan69}. In the sequel we shall proceed with the $\pm i$ constant for the sake of simplicity. 

This leads us directly to an example of a Salem set. 

\begin{lemma} \label{salemexample} Let $E=\{(x, {|x|}^2): x \in {\Bbb F}^{d-1}_q\}$. Then $E$ is a Salem set. \end{lemma} 

To prove the lemma, observe that $\# E=q^{d-1}$. Furthermore, 
$$ \widehat{E}(m,t)=q^{-d} \sum_{x \in {\Bbb F}^{d-1}_q} e^{\frac{2 \pi i (x \cdot m+t{|x|}^2)}{q}}.$$

Using (\ref{Gauss}) and (\ref{gauss}) we see that 
$$ |\widehat{E}(m,t)| \lesssim q^{-d} q^{\frac{d-1}{2}}, $$ and the lemma is proved. 
 
\section{Estimation of the finite field analog of the Mattila integral and proof of Theorem \ref{mattila} and Theorem \ref{salem}}

Using (\ref{distcoef}), (\ref{Gauss}), (\ref{gauss}), and (\ref{gaussd}), we see that if $k \not=0$, we have 
$$ \widehat{\nu}(k)=q^{-1} \frac{1}{{(\# E)}^2}  \sum_{x \in {\Bbb F}_q^d} T^k_qE(x) \cdot E(x)$$
$$=q^{-1} \frac{1}{{(\# E)}^2} \sum_{m \in {\Bbb F}^d_q} \sum_{x \in {\Bbb F}^d_q} e^{\frac{2 \pi i x \cdot m}{q}} \widehat{T^k_qE}(m) \cdot E(x)$$
$$=q^{-1} \frac{1}{{(\# E)}^2}  \sum_{m \in {\Bbb F}^d_q} \sum_{x \in {\Bbb F}^d_q} e^{\frac{2 \pi i x \cdot m}{q}} G(m,k)\widehat{E}(m) \cdot E(x)$$
$$=q^{-1} \frac{1}{{(\# E)}^2} \sum_{m \in {\Bbb F}^d_q} \sum_{x \in {\Bbb F}^d_q} e^{\frac{2 \pi i x \cdot m}{q}} E(x) \cdot \widehat{E}(m) e^{\frac{2 \pi i {|m|}^2}{4kq}} g^d(k)$$
$$ =q^{-1} \frac{1}{{(\# E)}^2}  q^d \sum_{m \in {\Bbb F}^d_q} {|\widehat{E}(m)|}^2 e^{\frac{2 \pi i {|m|}^2}{4kq}} g^d(k)$$
$$=q^{-1} q^{\frac{d}{2}} q^d \frac{1}{{(\# E)}^2}  {(\pm i)}^d \sum_{m \in {\Bbb F}^d_q} {|\widehat{E}(m)|}^2 e^{\frac{2 \pi i {|m|}^2}{4kq}} .$$

Squaring both sides, we get
$$ {|\widehat{\nu}(k)|}^2=q^{d-2} q^{2d} \frac{1}{{(\# E)}^4}  \sum_{m,m' \in {\Bbb F}^d_q} {|\widehat{E}(m)|}^2 {|\widehat{E}(m')|}^2 e^{\frac{2 \pi i( {|m|}^2-{|m'|}^2)}{4kq}}.$$

We conclude that 
$$ \sum_{k \in {\Bbb F}_q} {|\widehat{\nu}(k)|}^2=q^{-2}+q^{d-2} q^{2d} \frac{1}{{(\# E)}^4}  \sum_{k \in {\Bbb F}_q \backslash 0} \sum_{m,m' \in {\Bbb F}^d_q} {|\widehat{E}(m)|}^2 {|\widehat{E}(m')|}^2 e^{\frac{2 \pi i k( {|m|}^2-{|m'|}^2)}{q}}=$$
$$=\frac{q^{3d-1}}{{(\# E)}^4} \sum_{\{(m,m') \in {\Bbb F}^d_q \times {\Bbb F}^d_q: {|m|}^2={|m'|}^2 \}} {|\widehat{E}(m)|}^2 {|\widehat{E}(m')|}^2$$
$$+\frac{q^{3d-2}}{{(\# E)}^4} \sum_{m,m' \in {\Bbb F}^d_q} {|\widehat{E}(m)|}^2 {|\widehat{E}(m')|}^2.$$ 

Now, 
$$ \sum_{m,m' \in {\Bbb F}^d_q} {|\widehat{E}(m)|}^2 {|\widehat{E}(m')|}^2$$
$$={\left( \sum_{m \in {\Bbb F}^d_q} {|\widehat{E}(m)|}^2 \right)}^2$$
$$=q^{-2d} {\left( \sum_{x \in {\Bbb F}^d_q} E^2(x) \right)}^2 \approx q^{-2d} {(\# E)}^2.$$ 

We conclude that 
$$ \sum_{k \in {\Bbb F}_q} {|\widehat{\nu}(k)|}^2=\frac{q^{3d-1}}{{(\# E)}^4}  \sum_{\{(m,m') \in {\Bbb F}^d_q \times {\Bbb F}^d_q: {|m|}^2={|m'|}^2 \}} {|\widehat{E}(m)|}^2 {|\widehat{E}(m')|}^2+O \left(q^{-2} \frac{q^{d}}{{(\# E)}^2}\right), $$ and Theorem \ref{mattila} follows since 
$$ {\cal M}(q)=q^2 \sum_{k \in {\Bbb F}_q} {|\widehat{\nu}(k)|}^2. $$

To prove Theorem \ref{salem} observe that if the estimate (\ref{salemest}) holds, then  
$$ {\cal M}(q) \lesssim \frac{q^{3d+1}}{{(\# E)}^4}  \# \{(m,m') \in {\Bbb F}^d_q \times {\Bbb F}^d_q: {|m|}^2={|m'|}^2\} \cdot q^{-4d} {(\# E)}^2$$ 
$$=q^{d} {(\# E)}^{-2} \lesssim  1$$ if $\# E \gtrsim q^{\frac{d}{2}}$. This completes the proof. 

\vskip.125in 

\section{The finite field analog of the spherical average-proof of Theorem \ref{trivial}, \ref{plug}, and \ref{something} } 

\subsection{Proof of Theorem \ref{something} and Theorem \ref{trivial}} We have 
$$ \sigma^2_E(m)=\sum_{\{m' \in {\Bbb F}^d_q: {|m'|}^2={|m|}^2\}} {|\widehat{E}(m')|}^2$$
$$=q^{-2d} \sum_{x,y \in {\Bbb F}^d_q} E(x)E(y) \sum_{\{m' \in {\Bbb F}^d_q: {|m'|}^2={|m|}^2\}} e^{\frac{2 \pi i (x-y) \cdot m'}{q}}$$
$$=q^{-d} \sum_{x,y \in {\Bbb F}^d_q} E(x)E(y) \widehat{S}_{{|m|}^2}(x-y)$$
\begin{equation} \label{keyav} \lesssim q^{-\frac{d+1}{2}} \frac{{(\# E)}^2}{q^d}, \end{equation} which proves Theorem \ref{something} provided we can establish the following estimate. 

\begin{lemma} \label{salemsphere} The sphere $S_r$, $r \not=0$, is a Salem set. In other words, for any non-zero $x \in {\Bbb F}^d_q$, 
$$ |\widehat{S}_r(m)| \lesssim q^{-d} q^{\frac{d-1}{2}}, $$ and 
$$ \# S_r \approx q^{d-1}.$$
\end{lemma}  

Plugging (\ref{keyav}) into (\ref{mattilaintegral}) we get 
$$ \frac{q^{3d+1}}{{(\# E)}^4} q^{-\frac{d+1}{2}} \frac{{(\# E)}^2}{q^d} \sum_{m \in {\Bbb F}^d_q} {|\widehat{E}(m)|}^2$$
$$=\frac{q^{3d+1}}{{(\# E)}^2} q^{-\frac{d+1}{2}} q^{-d} \sum_{x \in {\Bbb F}^d_q} E^2(x)$$
$$=\frac{q^{2d+1} q^{-d} q^{-\frac{d+1}{2}}}{\# E}=\frac{q \cdot q^{\frac{d-1}{2}}}{\# E}.$$

This implies that 
$$ \# \Delta(E) \gtrsim \frac{\# E}{q^{\frac{d-1}{2}}}, $$ thus establishing Theorem \ref{trivial} up to the proof of Lemma \ref{salemsphere}. 

\vskip.125in

To prove Lemma \ref{salemsphere}, we write  
$$ \widehat{S}_r(m)=q^{-d} \sum_{\{x \in {\Bbb F}^d_q: {|x|}^2=r\}} e^{-\frac{2 \pi i x \cdot m}{q}}$$
$$=q^{-d} \sum_{x \in {\Bbb F}^d_q} q^{-1} \sum_{j \in {\Bbb F}_q} e^{\frac{2 \pi i j({|x|}^2-r)}{q}} e^{-\frac{2 \pi i x \cdot m}{q}}$$
$$=q^{-d-1} \sum_{j \in {\Bbb F}^{*}_q} e^{-\frac{2 \pi i jr}{q}} \sum_{x \in {\Bbb F}^d_q}  e^{\frac{2 \pi i j{|x|}^2}{q}} e^{-\frac{2 \pi i x \cdot m}{q}}$$
$$=q^{-d-1}  \sum_{j \in {\Bbb F}^{*}_q} e^{-\frac{2 \pi i jr}{q}}  G(-m,-j)$$
$$=q^{-d-1}  \sum_{j \in {\Bbb F}^{*}_q} e^{-\frac{2 \pi i jr}{q}} {(\pm i)}^d q^{\frac{d}{2}} e^{-\frac{2 \pi i {|m|}^2}{4j}}$$
$$=q^{-\frac{d}{2}} q^{-1} {(\pm i)}^d \sum_{j \in {\Bbb F}^{*}_q} e^{-\frac{2 \pi i}{q}(jr+\frac{{|m|}^2}{4j})}.$$ 

This reduces the proof of Lemma \ref{salemsphere} to the following Kloosterman sum estimate due to Andre Weil (\cite{We48}). See, for example, \cite{IK04} for a nice proof. 

\begin{lemma} \label{kloosterman} If $q$ is a prime, then 
$$ \left| \sum_{j \in {\Bbb F}^{*}_q} e^{-\frac{2 \pi i}{q} (jr+j^{-1}r')} \right| \lesssim \sqrt{q}$$ for any $r,r' \in {\Bbb F}_q$. \end{lemma} 

We now prove that $\# S_r \approx q^{d-1}$. By above, 
$$ \sum_{x,y \in {\Bbb F}^d_q} {|\widehat{S}_r(x)|}^2=q^{-d} q^{-2} \sum_{x \in {\Bbb F}^d_q} \sum_{u,v \in {\Bbb F}^{*}_q} e^{\frac{2 \pi i}{q}(r(u-v)+{|x|}^2(u^{-1}-v^{-1})} $$
$$=q^{-d-2} \sum_{\{(u,v) \in {\Bbb F}^{*}_q \times {\Bbb F}^{*}_q: u \not=v\}} e^{\frac{2 \pi i (u-v)r}{q}} q^{\frac{d}{2}}$$
$$+q^{-2} \sum_{u \in {\Bbb F}^{*}_q} 1=O(q^{-1}).$$

It follows that 
$$ \# S_r=\sum_{y \in {\Bbb F}^d_q} S^2_r(x)=q^d \sum_{x \in {\Bbb F}^d_q} {|\widehat{S}_r(x)|}^2=O(q^{d-1}), $$ as desired. 

\subsection{Proof of Theorem \ref{plug}} We now prove Theorem \ref{plug}. If the estimate \ref{disavest} holds, then 
$$ \frac{q^{3d+1}}{{(\# E)}^4}  \sum_{\{(m.m') \in {\Bbb F}^d_q \times {\Bbb F}^d_q: {|m'|}^2={|m|}^2\}} {|\widehat{E}(m)|}^2 {|\widehat{E}(m')|}^2$$
$$ \lesssim \frac{q^{3d+1}}{{(\# E)}^4}  q^{-\beta} \sum_{m \in {\Bbb F}^d_q} {|\widehat{E}(m)|}^2$$
$$=\frac{q^{3d+1}}{{(\# E)}^4} q^{-\beta} q^{-d} \sum_{x \in {\Bbb F}^d_q} E^2(x)=\frac{q^{2d+1-\beta}}{{(\# E)}^3}, $$ and the proof is complete. 

\vskip.125in 

\section{Non-Salem sets and the behavior of the Mattila integral} 

We have already seen that if $E$ is Salem set, then ${\cal M}(q) \lesssim 1$. Thus it makes sense to look at the behavior of ${\cal M}(q)$ in the case when $E$ is a not a Salem set. Let 
$$ E=\{(k,k): k \in {\Bbb F}_q\}.$$ 

It is immediately apparent that $\# \Delta(E)=q$. However, the point here is to show that $E$ is not a Salem set and that ${\cal M}(q)$ is nevertheless bounded. We have 
$$ \widehat{E}(m)=q^{-2} \sum_{k \in {\Bbb F}_q} e^{-\frac{2 \pi i (m_1+m_2)k}{q}}$$
$$=q^{-1} E'(m),$$ where 
$$ E'=\{(t,-t): t \in {\Bbb F}_q\}.$$ 

This shows that $E'$ is not a Salem set. On the other hand,
$$ {\cal M}(q)=q^{-1} \sum_{\{(m,m') \in {\Bbb F}^2_q \times {\Bbb F}^2_q: {|m|}^2={|m'|}^2\}} E'(m) \cdot E'(m')$$
$$=q^{-1} \sum_{\{(u,-u,v,-v): u,v \in {\Bbb F}_q, u^2=v^2\}} 1 \lesssim 1. $$

This example easily generalizes to higher (even) dimension. It is worth mentioning that this example is quite analogous to the Fourier transform of the Lebesgue measure on the boundary of a polygon in the plane. The Fourier transform behaves badly in directions normal to the sides of the polygon, but the decay rate is excellent away from those directions. The Mattila integral measures "average" decay of the Fourier transform and it is reasonable to conjecture that ${\cal M}(q)$ is bounded for all sets of cardinality $\gtrsim q^{\frac{d}{2}}$.

\end{document}